\documentclass{article}

\usepackage{graphicx}

\usepackage{amssymb,subfigure}
\usepackage{amsmath}

\newtheorem{teo}{Theorem}

\newtheorem{df}{Definition}
\newtheorem{rem}{Remark}
\newtheorem{lem}{Lemma}

\begin{document}

\title{Periodic wavelet frames and time-frequency localization
\footnote{The first author is supported by grant of President RF MK-1847.2012.1, by RFBR 12-01-00216-a, 
and by DAAD scholarship A/08/79920.}}

\author{Elena A. Lebedeva, J\"urgen Prestin}

\date{ealebedeva2004@gmail.com, prestin@math.uni-luebeck.de}

\maketitle

%\small

\begin{abstract}
A family of Parseval  periodic wavelet frames is constructed. The family 
has optimal time-frequency localization (in the sense of the Breitenberger uncertainty constant) 
with respect to a family parameter and it has the best currently known localization with respect 
to a multiresolution analysis parameter.
\end{abstract}



\section{Introduction}
\label{intr}

In recent years the wavelet theory of periodic functions has been continuously refined.
First, periodic wavelets were generated by periodization of wavelet functions on the real line 
(see, for example, \cite{da92}). A wider and more natural approach providing a flexibility on a 
theoretical front and in applications is to study periodic wavelets directly using a periodic analog 
of a multiresolution analysis (MRA). The concept of periodic 
MRA is introduced and discussed in \cite{klt95, MS04, NPS, P, PT, S98, S97,Zh}.
In \cite{GT1}, a unitary extension principle (UEP) for constructing Parseval wavelet frames 
is rewritten for periodic functions (see Theorem \ref{UEP}). The approach is developed further in \cite{GHS}. 
 
In this paper we focus on a property of good localization of both periodic wavelet 
functions and their Fourier coefficients. 
The quantitative characteristic of this property is an uncertainty constant (UC).
Originally, the concept of the UC was introduced for the real line case in 1927 
(see Definition \ref{Huc}) by Heisenberg in \cite{H}. 
Its periodic counterpart was introduced in 1985 by Brei\-tenberger in \cite{B}   (see Definition \ref{uc}).      
The smaller UC corresponds to the better localization. In both cases there exists 
a universal lower bound for the  UC (the  uncertainty principle). 
In non-periodic setup the minimum is attained on the Gaussian function. 
But there is no periodic function attending the lower bound.
 So, to find a sequence of periodic  functions having an 
asymptotically minimal UC and some additional setup, for example a wavelet structure, is a natural concern.
    
There is a connection between the Heisenberg and the Brei\-tenberger UCs for wavelets. 
In \cite{prqurase03} it is proved  that for periodic wavelets generated
by periodization (see the definition in  Section \ref{disc}) of a wavelet 
function on the real line the periodic UC tends
to the real line UC of the original function as a parameter of periodization tends to infinity. 
It would be a possible way to construct an optimal periodic wavelet system using the periodization 
of a wavelet system on the real line. 
However,  in \cite{Bat} and \cite{Balan} the following result is proven:
If a real line function $\psi$ generates a wavelet Bessel set and the frequency center 
$\omega_{0,\widehat{\psi^0}}=(\psi',\,\psi)_{L_2(\mathbb{R})}=0$ 
(see notation  $\omega_{0,\widehat{\psi^0}}$ in Definition \ref{Huc}), then the Heisenberg
UC is greater or equal to $3/2.$
Moreover, it is unknown if there exists a real line orthonormal wavelet 
basis or tight frame possessing the Heisenberg UC less than $2.134.$ 
This value is attained for a Daubechies wavelet \cite{GB}.  
The smallest possible  value of the Heisenberg UC for
the family  of the Meyer wavelets equals to $6.874$ \cite{LP}. 
It is well known \cite{ChW} that the Heisenberg UC of the  Battle-Lemarie and 
the Daubechies wavelets tends to infinity as their orders grow. A set of real line orthogonal 
wavelet bases with the uniformly bounded Heisenberg UCs as their orders (smoothness) 
grow is constructed in \cite{L08, L11}. 
On the other hand, there are examples of real line   
wavelet frames possessing asymptotically optimal UC such as
nonorthogonal B-spline wavelets \cite{UAE} and their generalizations \cite{GL}. 
However, these frames  are not tight and we are looking for an orthogonal basis or tight frame. 
We will discuss a particular issue of periodization in Section \ref{disc}.    
   
Some papers dealing with periodic UCs directly include \cite{GY, PQ, R05, Se}.
For the first time in \cite{Se} periodic UCs  uniformly bounded with respect to an MRA parameter 
are computed for so-called trigonometric wavelets (see also \cite{R05}). 
In \cite{GY}, it is shown that the UCs of uniformly local, regular, and stable periodic scaling 
functions and wavelets are uniformly bounded.
In \cite{PQ} an example of an asymptotically optimal set of periodic functions $\{\varphi_h\}_{h>0}$ 
is constructed, namely $UC(\varphi_h) < 1/2+\sqrt{h}/2.$
Later,  $\varphi_h$ is used as a scaling function to generate a 
stationary interpolatory  MRA $(V_n)$. For the corresponding wavelet functions $\psi_{n,h}$ the 
UC is optimal for a fixed space $V_n,$ but the estimate is nonuniform with respect to $n,$ namely 
$UC(\psi_{n,h})< 1/2+1.1 n^2 \sqrt{h}.$
Nothing changes after orthogonalization: $UC(\psi^{\bot}_{n,h})< 1/2+1.1 n^2 \sqrt{h},$\quad
$UC(\varphi^{\bot}_{n,h})< 1/2+ n^2 \sqrt{h}.$

The main contribution of this paper is Theorem \ref{opt}, where we construct a family of scaling sequences 
$\Phi^0=\left\{(\varphi_{j}^a)_j\ :\ a>1\right\}$ generating  a family of wavelet sequences
$\Psi^0=\left\{(\psi_{j}^a)_j\ :\ a>1\right\}$ corresponding to a nonstationary periodic MRA as 
it is defined in \cite{GHS}, \cite{klt95}, and \cite{S97}.  
For a fixed level $j$ of the MRA ($V_{2^j}$), similar to the construction in  \cite{PQ}, 
the UCs of $\varphi_{j}^a$ and $\psi_{j}^a$ are asymptotically optimal, that is
$$
\lim_{a \to \infty} \sup_{j \in \mathbb{N}} UC(\varphi_{j}^a)=\frac12,
\quad
\lim_{a \to \infty} UC(\psi_{j}^a)=\frac12.
$$ 
But now, for a fixed value of the parameter $a>1,$  the scaling sequence has the asymptotically
optimal UC, and the wavelet sequence has the smallest currently known value of 
the UC for the periodic wavelet frames setup, that is    
$$
\lim_{j \to \infty} \sup_{a>1} UC(\varphi_{j}^a)=\frac12,
\quad
\lim_{j \to \infty} UC(\psi_{j}^a)=\frac32.
$$
As it is indicated above, the functions constructed in \cite{PQ} do not have this property.

This issue partly answers the question stated in \cite{PQ} whether there exists a translation-invariant 
basis of a wavelet space  $W_j$ which is  asymptotically optimal independent of the MRA level $j.$ 
In Theorem \ref{opt} we get an 
affirmative answer for the case of scaling functions corresponding to tight wavelet frames. 
The case of wavelet basis is an open problem and it is a task for future work.
In this direction, some useful properties of shifted
Gaussian are discussed in \cite{KMNS}.  
We will consider the particular issue of wavelet sequences in Section \ref{disc}.

\section{Notations and auxiliary results}
\label{note}
Let $L_2(0,\, 1)$ be the space of all $1$-periodic square-integrable complex-valued functions, with
inner product $(\cdot,\cdot)$ given by 
$
(f,\,g):=\int_0^1 f(x)\overline{g(x)}\,\mathrm{d}x
$
for any $f,g, \in L_2(0,\,1),$ and norm $\|\cdot\|:=\sqrt{(\cdot,\,\cdot)}.$ 
The Fourier series of a function 
$
f \in L_2(0,\,1)
$
is defined by
$\sum_{k \in \mathbb{Z}}\widehat{f}(k) \mathrm{e}^{2 \pi \mathrm{i} k x},$
where its Fourier coefficient is defined by
$
\widehat{f}(k) = \int_0^1 f(x)\mathrm{e}^{-2 \pi \mathrm{i} k x}\,\mathrm{d}x.
$

Let $H$ be a separable Hilbert space. If there exist  constants $A,\,B>0$
such that for any $f \in H$ the following inequality holds
$
A \|f\|^2 \leq \sum_{n=1}^{\infty} \left|(f,\,f_n)\right|^2 \leq B \|f\|^2,
$ 
then the sequence $(f_n)_{n \in \mathbb{N}}$ is called  \texttt{a frame}  for $H.$
If $A=B(=1),$ then the sequence $(f_n)_{n \in \mathbb{N}}$ is called \texttt{a tight frame (a Parseval frame)} for $H.$ 
In addition, if $\|f_n\|=1$ for all $n \in \mathbb{N}$,  then the system forms an orthonormal basis.
More information about frames can be found in \cite{ch03}. 

In the sequel,  we use the following notation 
$
f_{j,k}(x):=f_j(x- 2^{-j} k)
$
for a function $f_j \in L_2(0,\,1).$
Consider  functions   $\varphi_0,\,\psi_j \in L_2(0,\,1),$ $j=0,1,\dots$ If the collection 
 $\Psi:=\left\{\varphi_0,  \psi_{j,k} : \ j=0,1,\dots,\ k=0,\dots,2^j-1 \right\},$
forms a frame (or basis) for $L_2(0,\,1)$ then $\Psi$ is said to be \texttt{a periodic wavelet frame 
(or wavelet basis)} for $L_2(0,\,1).$
Let us recall the UEP for a periodic setting. We consider a case of one wavelet generator.

\begin{teo}[\cite{GT1}]
\label{UEP}
Let $\varphi_j \in L_2(0,\,1),$ $j=0,1,\dots$,
be  a sequence of $1$-periodic functions such that 
\begin{equation}\label{con1}
     \lim_{j \to \infty}2^{j/2} \widehat{\varphi}_j(k) = 1.
\end{equation}
Let $\mu^j_k$ be a two-parameter sequence    such that $\mu^j_{k+2^j}=\mu^j_{k},$
and
\begin{equation}
\label{con2}
	\widehat{\varphi}_{j-1}(k)=\mu^j_k \widehat{\varphi}_{j}(k). 
\end{equation}
Let $\psi_j,$ $j=0,1,\dots,$ 
be a sequence of $1$-periodic functions defined using Fourier coefficients
\begin{equation}
\label{con3}
		{\widehat{\psi}}_{j}(k)=\lambda^{j+1}_{k} \widehat{\varphi}_{j+1}(k),
\end{equation}
where $\lambda^j_{k+2^j}=\lambda^j_{k}$  and 
\begin{equation}
\label{con4}
	\left(
	\begin{array}{cc}
	\mu ^j_k &  \mu ^j_{k+2^{j-1}} \\
	 \lambda^j_{k} & \lambda^j_{k+2^{j-1}}
	\end{array}
	\right)
	\left(
	\begin{array}{cc}
	\overline{\mu} ^j_k &  \overline{\lambda}^j_{k} \\
	\overline{\mu} ^j_{k+2^{j-1}} & \overline{\lambda}^j_{k+2^{j-1}}
	\end{array}
	\right)
	=
	\left(
	\begin{array}{cc}
	2 & 0 \\
	0 & 2
	\end{array}
	\right).
\end{equation}
Then the family
$\Psi : =\left\{\varphi_0,  \psi_{j,k} : \ j=0,1,\dots,\ k=0,\dots,2^j-1\right\}$
forms a Parseval wavelet frame for $L_2(0,\,1).$ 
\end{teo}
The sequences $(\varphi_j)_j,$ $(\psi_j)_j,$ $(\mu^j_k)_k,$ and $(\lambda^j_k)_k$ are 
called \texttt{scaling sequence, wavelet sequence, mask and wave\-let mask}   respectively.
This setup generates a periodic MRA: By definition, put\\
$V_j= {\rm span} \left\{\varphi_{j,k}; k= 0,\ldots,2^j-1\right\}$ 
for $j \geq 0.$ Then the sequence $(V_j)_{j \geq 0}$ is a periodic MRA.  
 
Let us recall the definitions of the  UCs and the uncertainty principles.
 
 \begin{df}[\cite{H}]
 \label{Huc}
\texttt{The (Heisenberg) UC} of $f \in L_2(\mathbb{R})$ is the functional 
$UC_H(f):=\Delta_{f}\Delta_{\widehat{f}}$ such that
$$
\begin{array}{ll}
\Delta_{f}^2:=\|f\|^{-2}_{L^2(\mathbb{R})}\int_{\mathbb{R}}(t-t_{0f})^2|f(t)|^2\,\mathrm{d}t, &
\Delta_{\widehat{f}}^2:=\|\widehat{f}\|^{-2}_{L^2(\mathbb{R})}
 \int_{\mathbb{R}}(\omega-\omega_{0\widehat{f}})^2|\widehat{f}(\omega)|^2\,\mathrm{d}\omega, \\
\end{array}
$$
$$
\begin{array}{ll}
t_{0f}:=\|f\|^{-2}_{L^2(\mathbb{R})}\int_{\mathbb{R}}t|f(t)|^2\,\mathrm{d}t, &
\omega_{0\widehat{f}}:=\|\widehat{f}\|^{-2}_{L^2(\mathbb{R})}\int_{\mathbb{R}}
\omega|\widehat{f}(\omega)|^2\,\mathrm{d}\omega. \\
\end{array}
$$  
\end{df} 
\begin{teo}[\cite{H}]
Let $f \in L_2(\mathbb{R})$, then $UC_H(f)\geq 1/2,$ and the equality is attained iff $f$ is the Gaussian. 
\end{teo}

\begin{df}[\cite{B}]
\label{uc}
Let $f(x)=\sum_{k \in \mathbb{Z}} c_k \mathrm{e}^{2 \pi \mathrm{i} k x}\in L_2(0,\,1).$  
\texttt{ The first trigonometric moment} is defined as
$$ 
\tau(f):=-2 \pi \int_0^1 \mathrm{e}^{2 \pi \mathrm{i} x} |f(x)|^2\, \mathrm{d}x =
-2\pi \sum_{k \in \mathbb{Z}} c_k \bar{c}_{k+1}.
$$
\texttt{The angular variance} of the function $f$ is defined by  
$$
{\rm var_A }(f):=\frac{1}{4\pi^2}\left( \frac{\left(\sum_{k \in \mathbb{Z}}|c_k|^2\right)^2}{
\left|\sum_{k \in \mathbb{Z}}c_k \bar{c}_{k+1}\right|^2}-1\right)
=
\frac{\|f\|^4}{|\tau(f)|^2}-\frac{1}{4\pi^2}.
$$
\texttt{The frequency variance} of the function $f$ is defined by  
$$
{\rm var_F }(f):= \frac{4\pi^2\sum_{k \in \mathbb{Z}}k^2 |c_k|^2}{\sum_{k \in \mathbb{Z}}|c_k|^2}-
\frac{4\pi^2\left(\sum_{k \in \mathbb{Z}}k|c_k|^2\right)^2}{\left(\sum_{k \in \mathbb{Z}}|c_k|^2\right)^2}
=
\frac{\|f'\|^2}{\|f\|^2}+\frac{(f',\, f)^2}{\|f\|^4}.
$$
The quantity 
$
UC(\{c_k\}):=UC(f):=\sqrt{\mathrm{var_A}(f)\mathrm{var_F}(f)}
$
is called \texttt{the periodic (Breitenberger) UC}.
\end{df}

\begin{teo}[\cite{B, PQ}]
\label{UC}
Let $f \in L_2(0,\,1)$,   $f(x)\neq C 
\mathrm{e}^{2 \pi \mathrm{i} k x},$ $C \in \mathbb{R},$ $k \in \mathbb{Z}$. Then
 $UC(f) > 1/2$  and there is no function such that  $UC(f) = 1/2.$
\end{teo}

Since periodic wavelet bases and frames are nonstationary in nature and the UC 
has no extremal function, it is natural to give the following
\begin{df}
\label{dfopt}
Suppose that $\varphi_j$\  ($\psi_j$) is a scaling (a wavelet) sequence. Then
the quantity 
$$
\limsup_{j \to \infty}UC(\varphi_j) \quad (\limsup_{j \to \infty}UC(\psi_j))
$$
is called \texttt{the UC of the scaling (the wavelet) sequence}. We say that
a sequence of periodic functions $(f_j)_{j \in \mathbb{N}}$
has \texttt{an optimal UC} if  
$$
\lim_{j \to \infty}UC(f_j)=1/2.
$$ 
\end{df} 
To justify the definition we note that since $\inf UC(f)=1/2,$ it follows that if 
$\limsup_{j \to \infty}UC(f_j)=1/2,$ then $\lim_{j \to \infty}UC(f_j)=1/2.$ 
So in the optimal case one can use $\lim_{j \to \infty}$ instead of $\limsup_{j \to \infty}.$

\section{Main results}

In the following theorem we construct a family of periodic Parseval wavelet 
frames with the optimal UCs for scaling functions and currently the best known UCs for wavelets.

\begin{teo}\label{opt}
There exists  a family of periodic wavelet sequences
$\Psi_a:=\{(\psi_{j}^a)_j\}_a$  such that
for any fixed $a>1$ the system $\{\varphi_0^a\}\cup \{\psi_{j,k}^a:\ j=0,1,\dots, \ k=0,\dots,2^j-1\}$ 
forms a Parseval frame in $L_2(0,\,1)$  and
\begin{gather}
\label{juclim}
\lim_{j \to \infty} \sup_{a>1} UC(\varphi_{j}^a)=\frac12, \quad
\lim_{a \to \infty} \sup_{j \in \mathbb{N}} UC(\varphi_{j}^a)=\frac12,
\\
\label{auclim}
\lim_{j \to \infty} UC(\psi_{j}^a)=\frac32, \quad
\lim_{a \to \infty}  UC(\psi_{j}^a)=\frac12.
\end{gather}
Put $\varphi_0^a = 1.$
 Let $\nu^{j,a}_{k}$ be  a sequence given by $\nu^{1,a}_0=\nu^{1,a}_1=\sqrt{1/2}$ and
\begin{equation}
\label{nu^ja}
	\nu^{j,a}_{k}:=\left\{
	\begin{array}{cl}
	\mathrm{exp}\left(-\frac{k^2+a^2}{j(j-1)a}\right), & k=-2^{j-2}+1,\dots,2^{j-2},\\[2ex]
	\sqrt{1- \mathrm{exp}\left(-\frac{2((k-2^{j-1})^2+a^2)}{j(j-1)a}\right)}, 
      & k=2^{j-2}+1,\dots, 3 \times 2^{j-2},
	\end{array}
	\right.
\end{equation}
and extended $2^j$-periodic with respect to $k$.	Furthermore, we define 
$\widehat{\xi}_{j}^a(k):=\prod_{r=j+1}^{\infty}\nu^{r,a}_{k}.$ 
Then the scaling sequence, masks, wavelet masks and wavelet sequence are defined 
	 respectively as
\begin{gather}
        \widehat{\varphi_{j}^a}(k):=2^{-j/2}\widehat{\xi}_{j}^a(k),\qquad\quad 
	\mu^{j,a}_k:=\sqrt{2} \nu^{j,a}_k,\ 
	\notag\\
	\lambda^{j,a}_k:=e^{2\pi i 2^{-j}k}\mu^{j,a}_{k+2^{j-1}},\qquad \quad
	{\widehat{\psi}}_{j}^a(k):=\lambda^{j+1,a}_{k} \widehat{\varphi}_{j+1}^a(k).
\label{uep0}
\end{gather}
\end{teo}

\begin{rem}
\label{homo}
The UC  is a homogeneous functional, 
that is $UC(\alpha f) = UC(f)$ for  $\alpha \in \mathbb{R},$ so 
$UC(\varphi_{j}^{a}) = UC(2^{-j/2}\xi_{j}^a)=UC(\xi_{j}^a)$ and in the sequel we prove the equalities 
$\lim_{j \to \infty}\sup_{a>1}UC(\xi_{j}^a)=1/2$ and\\
$\lim_{a \to \infty}\sup_{j\in \mathbb{N}} UC(\xi_{j}^a)=1/2$ instead of (\ref{juclim}).
By analogy, let $\eta^a_j:=2^{j/2}\psi^a_j,$ then 
$UC(\psi^a_j)=UC(\eta^a_j).$
\end{rem}

To prove Theorem \ref{opt}, we need some technical Lemmas.

\begin{lem}
\label{continuous}
The UC is a continuous functional 
with respect to the norm $\|f\|_{W^2_1}:=\|f\|+\|f'\|.$
\end{lem}

\noindent
\textbf{Proof.}
Indeed, $\tau(f)$ and $(f',\,f)$ are continuous with respect to this norm.
Using  the Cauchy-Bunyakovskiy-Schwarz inequality, we immediately get
\begin{eqnarray*}
 \frac{1}{2\pi}\left|\tau(f)-\tau(g)\right| 
  \leq 
 \int_0^1  \left||f|^2-|g|^2\right|  
 =
   \left(\Bigl||f|-|g|\Bigr|,\,|f|+|g|\right)\leq \Bigl\||f|-|g|\Bigr\| \, 
   \Bigl\||f|+|g|\Bigr\|
   \leq
    \Bigl(\|f\|+\|g\|\Bigr) \|f-g\|_{W^2_1}; 
\\
\left|(f',f)-(g',g)\right| \leq  
|(f',f-g)| + |(f'-g',g)|
\leq
 \|f'\|\  \|f-g\| + \|f'-g'\|\  \|g\|
\leq
\max\Bigl\{\|f'\|,\|g\|\Bigr\} \|f-g\|_{W^2_1}.
\end{eqnarray*} 
It remains to note that the UC continuously depends on 
$\|f\|,$ $\|f'\|,$ $\tau(f),$ and $(f',\,f).$ \hfill{$\Box$}

\begin{lem}
\label{tech}
Suppose $\alpha, \beta, \gamma \in \mathbb{R},$ $m=0,1,\dots,$ and $\ 0<\,b\,<M,$ 
where $M$ is an absolute constant, then 
\begin{gather}
\sum_{k \in\mathbb{Z}}(\alpha k^2 + \beta k + \gamma)^m\, \mathrm{e}^{-b (\alpha k^2 +
 \beta k + \gamma)}= (-1)^m
\left(
\mathrm{exp}\left(-b\left(\gamma-\frac{\beta^2}{4\alpha}\right)\right)
\sqrt{\frac{\pi}{b \alpha}}
\right)^{(m)}_{b^m}
+\mathrm{exp}\left(-\frac{\pi^2-\varepsilon}{b\alpha}\right)\,O(1),
\label{estpoi}
\end{gather}
as $b\to 0,$
where $\varepsilon>0$ is an arbitrary small parameter.
\end{lem}

\noindent
\textbf{Proof.}
It is possible to change the order of summation and differentiation, so
$$
\sum_{k \in\mathbb{Z}}(\alpha k^2 + \beta k + \gamma)^m\, 
\mathrm{e}^{-b (\alpha k^2 + \beta k + \gamma)}=(-1)^m
\left(
\sum_{k \in\mathbb{Z}}
\mathrm{e}^{-b (\alpha k^2 + \beta k + \gamma)}
\right)^{(m)}_{b^m}
$$
$$=
(-1)^m
\left(
\mathrm{exp}\left(-b  \left(\gamma - \frac{\beta^2}{4\alpha}\right)\right)
\sum_{k \in\mathbb{Z}}
\mathrm{exp}\left(-b \alpha \left(k + \frac{\beta}{2\alpha}\right)^2\right)
\right)^{(m)}_{b^m}.
$$
Using the Poisson summation formula for the function $f(t)=\mathrm{e}^{-b \alpha t^2}$ 
\begin{equation}
	\label{poi}
	\sum_{k \in \mathbb{Z}} \mathrm{e}^{-b \alpha(k-t)^2}=
   \sqrt{\frac{\pi}{b \alpha}}\sum_{k \in \mathbb{Z}} \cos 2\pi k t\ 
     \mathrm{exp}\left(\frac{-\pi^2 k^2}{b \alpha}\right),
\end{equation}
with  $t=-\beta/(2\alpha),$
and then differentiating $m$ times with respect to $b$, we get  
\begin{eqnarray*}
\lefteqn{(-1)^m \left(
\mathrm{exp}\left(-b  \left(\gamma - \frac{\beta^2}{4\alpha}\right)\right)
\sqrt{\frac{\pi}{b \alpha}}
\sum_{k \in\mathbb{Z}}
\cos\left(2\pi k \frac{\beta}{2\alpha}\right)\mathrm{exp}\left(-\frac{\pi^2 k^2}{b \alpha}\right)
\right)^{(m)}_{b^m}} &&\\
 &=&
(-1)^m \left(
\mathrm{exp}\left(-b\left(\gamma-\frac{\beta^2}{4\alpha}\right)\right) 
\sqrt{\frac{\pi}{b \alpha}}
\right)^{(m)}_{b^m} +(-1)^m 2\left(
\mathrm{exp}\left(-b\left(\gamma-\frac{\beta^2}{4\alpha}\right)\right) 
\sqrt{\frac{\pi}{b \alpha}}
\right)^{(m)}_{b^m}
\sum_{k=1}^{\infty}
\cos\left(2\pi k \frac{\beta}{2\alpha}\right)\mathrm{exp}\left(-\frac{\pi^2 k^2}{b \alpha}\right) \\
&& + (-1)^m  \sum_{r=1}^{m} \binom{m}{r} \left(
\mathrm{exp}\left(-b\left(\gamma-\frac{\beta^2}{4\alpha}\right)\right) 
\sqrt{\frac{\pi}{b \alpha}}
\right)^{(m-r)}_{b^{m-r}}
\left(1+2\sum_{k=1}^{\infty}
\cos\left(2\pi k \frac{\beta}{2\alpha}\right)\mathrm{exp}\left(-\frac{\pi^2 k^2}{b \alpha}\right)
\right)^{(r)}_{b^r}.
\end{eqnarray*}

For $r=0,\dots,m,$ we estimate
\begin{eqnarray*}
\left|\left(
\sum_{k=1}^{\infty}
\cos\left(2\pi k \frac{\beta}{2\alpha}\right)\mathrm{e}^{-\pi^2 k^2/(b \alpha)}
\right)^{(r)}_{b^r}\right|\leq
\sum_{k=1}^{\infty}
Q_r\left(k,\frac{1}{b}\right)\mathrm{e}^{-\pi^2 k^2/(b \alpha)} &\leq&
\mathrm{e}^{- \pi^2 /(b \alpha)}
\sum_{k=1}^{\infty}
Q_r\left(k,\frac{1}{b}\right)\mathrm{exp}\left(-\frac{\pi^2 (k^2-1)}{M \alpha}\right)\\
 &=&
\mathrm{e}^{-(\pi^2-\varepsilon)/(b \alpha)} O(1),
\end{eqnarray*}
where
$Q_r(k,1/b)$ is a polynomial of degree $2r$  in $k,$ and $1/b.$ We estimate summands of the form
$\mathrm{e}^{-\pi^2/(b \alpha)}/b^{\xi},$ $0<\xi<2m$ by
$\mathrm{e}^{-\pi^2/(b \alpha)}/b^{\xi}< \mathrm{exp}\left(-\frac{\pi^2-\varepsilon}{b \alpha}\right).$  
 \hfill{$\Box$}

\begin{lem}
\label{expw}
Suppose 
$\eta^{a,0}_{j}(t):=\sum_{k \in \mathbb{Z}}\widehat{\eta}_{j}^{a,0} (k) 
\mathrm{e}^{2 \pi \mathrm{i} k t},$
where 
\begin{equation}
\label{eta}
 \widehat{\eta}_{j}^{a,0} (k):=\mathrm{e}^{2\pi \mathrm{i} 2^{-j-1}k} 
\sqrt{1-\mathrm{exp}\left(-\frac{2(k^2+a^2)}{(j(j+1)a)}\right)}\ \ 
\mathrm{exp}\left(-\frac{k^2+a^2}{(j+1)a}\right);
\end{equation}
then
$
\lim_{j \to \infty}UC(\eta_{j}^{a,0})=3/2
$ for any fixed $a>1$ and 
$
\lim_{a \to \infty}UC(\eta_{j}^{a,0})=1/2
$
for any fixed $j \in \mathbb{N}.$  
\end{lem}

\noindent
\textbf{Proof.}
We estimate the quantities $((\eta_{j}^{a,0})',\,\eta_{j}^{a,0}),$ 
$\|\eta_{j}^{a,0}\|^2,$ $\|(\eta_{j}^{a,0})'\|^2,$ and 
$|\tau(\eta_{j}^{a,0})|$ and then substitute the expressions in Definition \ref{uc}.
Since $|\widehat{\eta}_{j}^{a,0} (k)|=|\widehat{\eta}_{j}^{a,0} (-k)|,$ we see that 
$
((\eta_{j}^{a,0})',\,\eta_{j}^{a,0})=\sum_{k \in \mathbb{Z}} k\ |\widehat{\eta}_{j}^{a,0} (k)|^2=0.
$

For convenience we replace $j+1$ by $1/h$  and $a$ by $1/q$. 
Then, $0<h \leq 1/2$, $0<q\leq 1$, $h\to 0,$ and $q\to 0.$ 
However, to avoid the fussiness of notations we keep 
the former name for the  function $\eta_{j}^{a,0}.$  By  (\ref{eta}),
\begin{eqnarray*}
\|\eta_{j}^{a,0}\|^2 =
\sum_{k \in \mathbb{Z}} |\widehat{\eta}_{j}^{a,0} (k)|^2  
=
\mathrm{exp}\left(-\frac{2h}{q}\right)\sum_{k \in \mathbb{Z}} 
\mathrm{exp}\left(-2 h q k^2\right)-
\mathrm{exp}\left(-\frac{2h}{(1-h)q}\right)\sum_{k \in \mathbb{Z}} 
\mathrm{exp}\left(-\frac{2 h q}{1-h} k^2\right).
\end{eqnarray*}
Using (\ref{estpoi}) twice  for $\alpha=1,$ $\beta=0,$ $\gamma=0,$ $m=0$,  
$b=2hq$ and $b=2hq/(1-h),$ we get
\begin{equation}
	\|\eta_{j}^{a,0}\|^2 =
  \mathrm{exp}\left(-\frac{2h}{q}\right)\sqrt{\frac{\pi}{2hq}}-
    \mathrm{exp}\left(-\frac{2h}{q(1-h)}\right) \sqrt{\frac{\pi(1-h)}{2hq}} +
    \left(\mathrm{e}^{C(h,\,q)}+\mathrm{e}^{C(h/(1-h),\,q)}\right)O(1),   
\label{norm}
\end{equation}
as $hq\to +0,$  where $C(h,\,q)=-2h/q-(\pi^2-\varepsilon)/(2hq).$

Similarly, to estimate the quantities $\|(\eta_{j}^{a,0})'\|^2,$ by (\ref{eta}),  
we write
\begin{equation*}
\frac{1}{4\pi^2}\|(\eta_{j}^{a,0})'\|^2 =
 \sum_{k \in \mathbb{Z}} k^2 |\widehat{\eta}_{j}^{a,0} (k)|^2 = 
  \mathrm{exp}\left(-\frac{2h}{q}\right)\sum_{k \in \mathbb{Z}} k^2 
  \mathrm{exp}\left(-2 h q k^2\right)
-\mathrm{exp}\left(-\frac{2h}{q(1-h)}\right)
\sum_{k \in \mathbb{Z}} k^2 \mathrm{exp}\left(-\frac{2 h q k^2}{1-h}\right).
\end{equation*}
Using (\ref{estpoi}) twice for $\alpha=1,$ $\beta=0,$ $\gamma=0,$ $m=1$,  
$b=2hq$ and $b=2hq/(1-h),$ we get
\begin{equation*} 
 \frac{1}{4\pi^2} \|(\eta_{j}^{a,0})'\|^2 
 =
 \frac{1}{2}
\mathrm{exp}\left(-\frac{2h}{q}\right)
	  \sqrt{\frac{\pi}{(2hq)^3}}-
	\frac{1}{2}
	 \mathrm{exp}\left(-\frac{2h}{q(1-h)}\right)
  \sqrt{\frac{\pi(1-h)^3}{(2hq)^3}}
 +\left( \mathrm{e}^{C(h,\,q)} + \mathrm{e}^{C(h/(1-h),\,q)}\right)O(1),
\end{equation*}
as $hq\to +0,$
where $C(h,\,q)$ is defined after formula (\ref{norm}).
So, recalling  $((\eta_{j}^{a,0})',\,\eta_{j}^{a,0})=0,$ by Definition \ref{uc}, 
we get the following asymptotic form for the frequency variance: 
\begin{equation}
      \label{vF}
	\frac{1}{4\pi^2}
	\frac{\|(\eta_{j}^{a,0})'\|^2}{\|\eta_{j}^{a,0}\|^2}\thicksim 
	\frac{3}{4hq}\ \ \mbox{ as } h \to 0
	\qquad \mbox{and}\qquad
	 \frac{1}{4\pi^2}
	\frac{\|(\eta_{j}^{a,0})'\|^2}{\|\eta_{j}^{a,0}\|^2}\thicksim 
	\frac{1}{4hq}\ \  \mbox{ as } q \to 0.
\end{equation}

To estimate the first trigonometric moment $\tau(\eta_{j}^{a,0})$ 
(see Definition \ref{uc}), by (\ref{eta}), we obtain
\begin{eqnarray*}
 \frac{1}{2\pi}|\tau(\eta_{j}^{a,0})|&=& 
\left|\sum_{k \in \mathbb{Z}}\widehat{\eta}_{j}^{a,0} (k) 
\overline{\widehat{\eta}_{j}^{a,0} (k+1)}\right|\\
  &=&
    \mathrm{e}^{-\frac{2h}{q}}\sum_{k \in \mathbb{Z}} 
    \sqrt{\left(1-\mathrm{exp}\left(-\frac{2h^2(k^2q^2+1)}{(1-h)q}\right)\right)
  \left(1-\mathrm{exp}\left(-\frac{2 h^2(q^2(k+1)^2+1)}{(1-h)q}\right)\right)}\
   \ \mathrm{e}^{-h q(2k^2+2k+1)} .
\end{eqnarray*}
Our task is to get the following representation for $|\tau(\eta_{j}^{a,0})|$:
\begin{eqnarray}
\label{tau_h!}
\frac{1}{2\pi}|\tau(\eta_{j}^{a,0})|
&=&\!\! \!\! \frac{\mathrm{e}^{-\frac{2h}{q}-\frac{hq}{2}}}{1-h} \sqrt{\frac{\pi}{8q}} 
 \left(\!\!\sqrt{h}+
\frac{(1-h)(16-4q^2)-3q}{4 q (1-h)}\sqrt{h^3}\right)+O(h^2 |\ln h|)\quad  
 \mbox{for a fixed } q\leq 1 \mbox{ and } h \to 0,
\\
\label{tau_q!}
\frac{1}{2\pi}|\tau(\eta_{j}^{a,0})|
&=&  \mathrm{e}^{-\frac{2h}{q}} \left(\mathrm{e}^{-\frac{h q}{2}}
\sqrt{\frac{\pi}{2 h q }}+O\Bigr(\frac{1}{\sqrt{q}}\mathrm{e}^{-\frac{2h^2}{q(1-h)}}\Bigr)\right)
\qquad\qquad\qquad\mbox{ for a fixed } h\leq 1/2 \mbox{ and } q \to 0.
\end{eqnarray}
Let us prove the estimate (\ref{tau_h!}).
Put
\begin{equation}
\label{notation}
	d:=\frac{2h^2}{1-h},\ \mathrm{v}(k):=q k^2+\frac{1}{q},\ s(k):=2k^2+2k+1.
\end{equation}
Thus, the first trigonometric moment is rewritten as follows:
\begin{equation}
	\label{tau0}
\frac{1}{2\pi}|\tau(\eta_{j}^{a,0})|=
\mathrm{e}^{-\frac{2h}{q}}\sum_{k \in \mathbb{Z}}
 \sqrt{(1-\mathrm{e}^{-d \mathrm{v}(k)})(1-\mathrm{e}^{-d \mathrm{v}(k+1)})}\mathrm{e}^{-h q s(k)}.
\end{equation}
Using the Taylor formula for the function 
$f(d)=\sqrt{(1-\mathrm{e}^{-d \mathrm{v}(k)})(1-\mathrm{e}^{-d \mathrm{v}(k+1)})}$ 
in the neighborhood of $d=0$, we get
$
f(d)=\sqrt{\mathrm{v}(k) \mathrm{v}(k+1)}d-\frac{1}{4}\sqrt{\mathrm{v}(k) 
\mathrm{v}(k+1)}(\mathrm{v}(k)+\mathrm{v}(k+1))d^2+ \frac{f'''(\bar{d})}{6} d^3,
$
and

\noindent
$
f'''(d)=
\frac{1}{8} N^{-\frac{5}{2}}\ M^{-\frac{5}{2}}
\bigg(\nu^3  N^3(1-M)(3+M^2)-
\mu \nu M N(1-M)(1-N)\Bigl(\mu M+\nu M+(\nu + \mu)M N\Bigr)+
\nu^3 M^3(1-N)(3+M^2)\bigg),
$
where $N:=1-\mathrm{e}^{-d \mathrm{v}(k)},$ $M:=1-\mathrm{e}^{-d \mathrm{v}(k+1)},$ 
$\nu:=\mathrm{v}(k),$ $\mu:=\mathrm{v}(k+1).$
We have 
$
|f'''(\bar{d})| d^3 = O(s^3(k) h^6).
$
Indeed, $f'''$ is a decreasing function on $0<d<1$. Collecting summands appropriately,  
one can check that $f''$ is a concave function on 
$0<d<1$. So, $|f'''(d)|\leq \lim_{d\to 0}f'''(d)=1/16 \sqrt{\mu \nu} (5 \mu^2 + 6 \mu \nu + 5 \nu^2).$ 
It remains to note that   
$\lim_{k\to \infty}\lim_{d\to 0}f'''(d)/s^3(k)=q^3/8$ is finite.

Applying (\ref{estpoi})  for $\alpha=1,$ $\beta=0,$ $\gamma=0,$ $m=3,$ $b=hq$, we have for the reminder of 
(\ref{tau0})
$$
\mathrm{e}^{-\frac{2h}{q}}\sum_{k \in \mathbb{Z}}f'''(\bar{d}) d^3 \mathrm{e}^{-h q s(k)} =
O\left(h^6\sum_{k \in \mathbb{Z}} s^3(k) \mathrm{e}^{-h q s(k)}\right)  
 = h^6  O\left(h^{-7/2} + \mathrm{e}^{-\frac{\pi^2-\varepsilon}{2h q}}O(1)\right)=O(h^{5/2}),
$$
as $h \to 0.$
Therefore, (\ref{tau0}) takes the form
$$
\frac{1}{2\pi}|\tau(\eta_{j}^{a,0})|
= \mathrm{e}^{-\frac{2h}{q}} 
\left| \sum_{k \in \mathbb{Z}} 
 \sqrt{\mathrm{v}(k) \mathrm{v}(k+1)} \left( d-\frac{1}{4}(\mathrm{v}(k)+\mathrm{v}(k+1))d^2 \right) 
  \mathrm{e}^{-h q s(k)}      \right|+O(h^{5/2}).
$$ 
With $u:=1/k$ we define the function $g$ by 
$$
g(u):=\frac{1}{k^2}
\sqrt{\mathrm{v}(k) \mathrm{v}(k+1)}=\frac{1}{k^2}
\sqrt{\left(q k^2+\frac{1}{q}\right)\left(q (k+1)^2+\frac{1}{q}\right)}
= \sqrt{\left(q+\frac{u^2}{q}\right)\left(q(u+1)^2+
\frac{u^2}{q}\right)}\ .
$$
Using the Taylor formula for  $g(u)$ in the neighborhood of $u=0$, we obtain
$
g(u)=q+q u +\frac{1}{q}u^2+\frac{g'''(\bar{u})}{6}u^3,
$
where
\begin{eqnarray*}\lefteqn{
g'''(u) = \frac{1}{2g(u)} \frac{6u (2/q+2q)+6(2u /q+2(1+u)q)}{q}
-\frac{3}{4g^3(u)}\left((u^2/q+q)(2/q+2q) \right.}&&\\
&+&\left.\frac{4u(2u/q+
2(1+u)q)+4(u^2/q+(1+u)^2q)}{q}\right)
\left((u^2/q+q)(2u/q+2(1+u)q)+\frac{2u(u^2/q+(1+u)^2 q)}{q}\right)\\
 &+&\frac{3}{8g^5(u)}\left((u^2/q+q)(2u/q+2(1+u)q)+
\frac{2u(u^2/q+(1+u)^2 q)}{q}\right)^3
=:\frac{1}{g(u)}P_1(u)+\frac{1}{g^3(u)}P_2(u)+\frac{1}{g^5(u)}P_3(u).
\end{eqnarray*}
Suppose $k \neq 0,$ then $-1 \leq u \leq 1.$
For a fixed $0<q\leq 1$, the value $|g'''(\bar{u})|$ is bounded.
Indeed, since $q+u^2/q \geq q$ and
$q(u+1)^2+u^2/q \geq q/(q^2+1),$
then
$0 < 1/g(u)\leq \sqrt{q^2+1}/q,$
 and the polynomials
$P_1(u),\ P_2(u),\ P_3(u)$ in $u$ are bounded on $[-1,\,1].$
So,  $g'''(\bar{u})u^3=u^3 O(1).$
Therefore,  we have
$$
\frac{1}{2\pi}|\tau(\eta_{j}^{a,0})|
=  \mathrm{e}^{-\frac{2h}{q}} 
\left| \sum_{k \in \mathbb{Z}, k \neq 0} 
 \left(k^2 q+ k q  +\frac{1}{q}+\frac{O(1)}{k}\right) 
 \left( d-\frac{1}{4}(\mathrm{v}(k)+\mathrm{v}(k+1))d^2 \right) 
  \mathrm{e}^{-h q s(k)}      \right|   +O(h^{2}).
$$
In the latter formula, we omit the summand for $k=0$ which equals to 
$
\mathrm{e}^{-\frac{2h}{q}} d \sqrt{1/q(q+1/q)}\left(1-
1/4 \ (2/q+q)d\right)\mathrm{e}^{-h q}=O(h^{2}).
$   
 Recalling (\ref{notation}), we estimate the coefficient of $O(1)$ in the  latter series 
\begin{eqnarray*}
A:=\left|\sum_{k =1}^{\infty}\frac{1}{k}\left( d-\frac{1}{4}(\mathrm{v}(k)+\mathrm{v}(k+1))d^2 \right) 
  \mathrm{e}^{-h q s(k)}  \right|
   \leq  d \sum_{k =1}^{\infty}\frac{1}{k} \mathrm{e}^{-h q s(k)}+\frac{q d^2}{4} 
   \sum_{k=1}^{\infty} \frac{s(k)}{k} \mathrm{e}^{- h q s(k)} +\frac{d^2}{2q} 
   \sum_{k=1}^{\infty}\frac{1}{k} \mathrm{e}^{- h q s(k)}.
 \end{eqnarray*}
 The first series is the main term as $h\to 0.$
Indeed, since $\sum_{k=1}^{\infty} \frac{s(k)}{k} \mathrm{e}^{- h q s(k)} < 
\sum_{k\in \mathbb{Z}}s(k) \mathrm{e}^{- h q s(k)}$ and 
 $\sum_{k=1}^{\infty}\frac{1}{k} \mathrm{e}^{- h q s(k)}<\sum_{k\in \mathbb{Z}} \mathrm{e}^{- h q s(k)},$ 
 then applying (\ref{estpoi}) for $\alpha=2,$ $\beta=2,$ $\gamma=1,$ $b=hq,$ $m=0,$ and $m=1$  we get 
 $d^2 \sum_{k\in \mathbb{Z}}s(k) \mathrm{e}^{- h q s(k)} \thicksim h^{5/2},$
 $
 d^2 \sum_{k\in \mathbb{Z}} \mathrm{e}^{- h q s(k)} \thicksim h^{7/2},
 $
 as $h\to 0.$
 Hence,
 \begin{eqnarray*}
 A   &\leq&  C_1 d \sum_{k =1}^{\infty}\frac{1}{k} \mathrm{e}^{-h q k^2}
   = C_1 d \left(\mathrm{e}^{-h q} + \sum_{k =2}^{\infty}\frac{1}{k} \mathrm{e}^{-h q k^2}\right)
   \leq C_1 d \left(\mathrm{e}^{-h q} + 
\int_{1}^{\infty}\frac{1}{x} \mathrm{e}^{-h q x^2}\, \mathrm{d}x\right)=
C_1 d \left(\mathrm{e}^{-h q} + \int_{\sqrt{h q}}^{\infty}\frac{1}{x} \mathrm{e}^{- x^2}\, 
\mathrm{d}x\right) \\
  &=&
 C_1 d \left(\mathrm{e}^{-h q} -\mathrm{e}^{-h^2 q^2} \ln (h q) + 
\int_{\sqrt{h q}}^{\infty}2 x \mathrm{e}^{- x^2} \ln x\, \mathrm{d}x\right)= O(h^2 |\ln h|).
\end{eqnarray*}
Similarly, one can estimate $\sum_{k < 0}$.
Finally, recalling (\ref{notation}), we have 
\begin{eqnarray}\nonumber
   \lefteqn{  \frac{1}{2\pi}|\tau(\eta_{j}^{a,0})|
	 =   \mathrm{e}^{-\frac{2h}{q}}} &&\\ \label{tau_h} \displaystyle
   \times &  \displaystyle	 \left|
	-\frac{d^2q^2}{8}\sum_{k \in \mathbb{Z}}  s^2(k)	 \mathrm{e}^{-h q s(k)}+
	 \frac{d}{8}(4q-4d+dq^2)\sum_{k \in \mathbb{Z}} s(k)	 \mathrm{e}^{-h q s(k)}+
	 \frac{d(2-q^2)(2q-d)}{4q^2}\sum_{k \in \mathbb{Z}}	 \mathrm{e}^{-h q s(k)}
	\right| +O(h^{2} |\ln h|).
\end{eqnarray}
Here, we return the summand for $k=0$, since it equals to
$d/q\ (1-(q+1/q)d/4)\mathrm{e}^{-hq}=O(h^2).$
To obtain (\ref{tau_h!}), it remains to substitute (\ref{estpoi})  for $\alpha=2,$ $\beta=2,$ 
$\gamma=1,$ $b=hq,$ $m=0,\,1,\,2$ in (\ref{tau_h}).

Using the Taylor formula for squared of (\ref{norm}) and (\ref{tau_h!}), we get
 \begin{eqnarray*}
 \frac{1}{4\pi^2}|\tau(\eta_{j}^{a,0})|^2=\frac{\pi}{8 q} h\left(1+\frac{8+q-6q^2}{2 q} h
 +O(h^{3/2}|\ln h|)\right) 
 \quad\mbox{and}\quad
 \|\eta_{j}^{a,0}\|^4=\frac{\pi}{8 q} h\left(1+\frac{8+q}{2 q} h
 +O(h^2)\right).
 \end{eqnarray*}
 Finally, substituting the last expressions and 
(\ref{vF}) in Definition \ref{uc} and calculating the limit we get
$
\lim_{j \to \infty} UC(\eta^{a,0}_j)=\frac32.
$ 
 
Let us prove (\ref{tau_q!}).
We start with (\ref{tau0}).
Using the mean value theorem for $f_0(x)=\sqrt{1-x},$ $x_0=0$, we obtain 
\begin{eqnarray*}
\lefteqn{
\frac{1}{2\pi}|\tau(\eta_{j}^{a,0})|=
\mathrm{e}^{-\frac{2h}{q}}\sum_{k \in \mathbb{Z}}\left(1-C_0(q,k) \mathrm{e}^{-d 
\mathrm{v}(k)}\right)\left(1-C_0(q,k+1) \mathrm{e}^{-d \mathrm{v}(k+1)}\right) \mathrm{e}^{-hq s(k)}
=
\mathrm{e}^{-\frac{2h}{q}}\left(\sum_{k \in \mathbb{Z}} \mathrm{e}^{-hq s(k)}\right.}&&\\
&\displaystyle - \left.
  \sum_{k \in \mathbb{Z}} C_0(q,k)\mathrm{e}^{-d \mathrm{v}(k)} \mathrm{e}^{-hq s(k)}-
 \sum_{k \in \mathbb{Z}}C_0(q,k+1)\mathrm{e}^{-d \mathrm{v}(k+1)} \mathrm{e}^{-hq s(k)}\right. \left.+
 \sum_{k \in \mathbb{Z}} C_0(q,k)C_0(q,k+1) 
 \mathrm{e}^{-d (\mathrm{v}(k)+\mathrm{v}(k+1))} \mathrm{e}^{-hq s(k)}\right),
\end{eqnarray*}
where $C_0(q,k):=1/(2\sqrt{1-c_0(q,k)}),$ $0<c_0(q,k)<\mathrm{e}^{-d \mathrm{v}(k)}.$
The first series is the main term as $q\to 0.$
Indeed, $C_0(q,k)$ is bounded 
 (for example, since $0<\mathrm{e}^{-d \mathrm{v}(k)}<1/2$ for a fixed $0<h\leq 1/2,$ and 
$0<q<(1-h)(2h^2)^{-1}\log 2,$ then $1/2 < C_0(q,k) < \sqrt{2}/2$).
So, the second, third, and fourth terms are estimated by 
$S_2:=\sqrt{2}/2\sum_{k \in \mathbb{Z}}\mathrm{e}^{-d \mathrm{v}(k)} \mathrm{e}^{-hq s(k)},$   
$S_3:=\sqrt{2}/2\sum_{k \in \mathbb{Z}}\mathrm{e}^{-d \mathrm{v}(k+1)} \mathrm{e}^{-hq s(k)},$
and
$S_4:=\sqrt{2}/2\sum_{k \in \mathbb{Z}}\mathrm{e}^{-d (\mathrm{v}(k)+\mathrm{v}(k+1))} \mathrm{e}^{-hq s(k)}$
respectively.   
Using (\ref{estpoi}) for an appropriate $\alpha,$ $\beta,$ $\gamma,$ $b,$ and $m=0,$ we see that 
$S_n=O(\frac{1}{\sqrt{q}}\mathrm{e}^{-\frac{d}{q}})$ for $n=2,3,4.$ 
To obtain (\ref{tau_q!}) it remains to apply  (\ref{estpoi}) (for $\alpha=2$, $\beta=2,$ 
$\gamma=1,$ $b=hq,$ $m=0$) to the first series  $\sum_{k \in \mathbb{Z}}\mathrm{e}^{-hq s(k)}.$

Finally, substituting 
(\ref{norm}), (\ref{vF}), and (\ref{tau_q!}) in Definition \ref{uc} and calculating the limit we obtain
$
\lim_{a \to \infty}  UC(\eta^{a,0}_j)=\frac12.
$
This completes the proof of Lemma \ref{expw}. \hfill$\Box$  
 
\medskip
 
\noindent
\textbf{Proof of Theorem \ref{opt}.} 
1. By Theorem \ref{UEP},  the family 
$\Psi_a : =\left\{ \textbf{1} ,  \psi_{j,k}^a : \ j=0,1,\dots,\ k=0,\dots,2^j-1\right\}$
(see (\ref{uep0}))
forms a Parseval wavelet frame for $L_2(0,\,1)$ for a  fixed $a>1.$
Indeed,
using definition (\ref{nu^ja}) and the elementary identity 
$
j^{-1}(j-1)^{-1}= (j-1)^{-1}-j^{-1},
$ 
we get 
\begin{equation}
	\label{xi0}
	\widehat{\xi}_{j}^a(k)
	=
	\left\{
	\begin{array}{ll}\displaystyle
	\prod_{r=j+1}^{J-1}\nu^{r,a}_{k} \prod_{r=J}^{\infty}\nu^{r,a}_{k}=
     \left(\prod_{r=j+1}^{J-1}\nu^{r,a}_{k}\right) \mathrm{exp}\left(-\frac{k^2+a^2}{(J-1)a}\right) , &
	j \leq J-2, \\[2ex] \displaystyle
	\prod_{r=j+1}^{\infty} \mathrm{exp}\left(-\frac{k^2+a^2}{r(r-1)a}\right)=
     \mathrm{exp}\left(-\frac{k^2+a^2}{ja}\right), &
	j > J-2, \\
	\end{array}
	\right.
\end{equation}
where $J=\lfloor\log_2(|k-1/2|+1/2)+3\rfloor.$ Therefore, the coefficients 
$\widehat{\xi}_{j}^a(k)$ are well-defined.
Then a straightforward calculation shows that conditions (\ref{con1})-(\ref{con4}) hold.
  
2.
According to Remark \ref{homo}, let us check $\lim_{j \to \infty}\sup_{a>1}UC(\xi_{j}^a)=1/2$ and
$\lim_{a \to \infty}\sup_{j\in \mathbb{N}} UC(\xi_{j}^a)=1/2$ instead of (\ref{juclim}).
Let us denote 
\begin{equation}
\label{ksi}
\xi_{j}^{a,0}(x):=\sum_{k \in \mathbb{Z}} \mathrm{e}^{-\frac{k^2+a^2}{ja}} 
\mathrm{e}^{2 \pi \mathrm{i} k x}=\mathrm{e}^{-\frac{a}{j}}
\sum_{k \in \mathbb{Z}} \mathrm{e}^{-\frac{k^2}{ja}} \mathrm{e}^{2 \pi \mathrm{i} k x}.
\end{equation}
Since the $UC$ is  homogeneous, it follows that 
$
UC(\xi_{j}^{a,0})=UC \left(\left\{\mathrm{e}^{-\frac{k^2}{j a}}\right\}\right).
$
It is known (see \cite{PQ}) that \\
$
\lim_{j \to \infty} UC\bigl(\bigl\{ \mathrm{e}^{-\frac{k^2}{j}}\bigr\} \bigr)=1/2.
$ 
Substituting $j a$ for $j$ to the last equality and swapping $j$ and $a$ we immediately get 
$$
\lim_{j \to \infty} \sup_{a>1} UC(\xi_{j}^{a,0})=\frac12,
\qquad
\lim_{a \to \infty} \sup_{j \in \mathbb{N}} UC(\xi_{j}^{a,0})=\frac12.
$$
So, taking into account the continuity of the UC (see Lemma \ref{continuous}), it remains to prove
that 
$\lim_{j \to \infty}\|\xi_{j}^{a}-\xi_{j}^{a,0}\|_{W^2_1}=0$ uniformly on $a>1,$
and 
$\lim_{a \to \infty}\|\xi_{j}^{a}-\xi_{j}^{a,0}\|_{W^2_1}=0$
uniformly on $j \in\mathbb{N}.$
Applying the elementary observation to  
$c_{j,k}=\widehat{\xi}_{j}^{a}(k)-\widehat{\xi}_{j}^{a,0}(k),$ (namely, 
if $\lim_{j \to \infty}c_{j,0}=0$ and
$\lim_{j \to \infty}\sum_{k \in \mathbb{Z}}k^2 |c_{j,k}|^2=0,$
then  $\lim_{j \to \infty}\sum_{k \in \mathbb{Z}} |c_{j,k}|^2=0$)
we see that it is sufficient to check  
$\lim_{j \to \infty}\|(\xi_{j}^{a})'-(\xi_{j}^{a,0})'\|=0,$ and 
$\lim_{a \to \infty}\|(\xi_{j}^{a})'-(\xi_{j}^{a,0})'\|=0.$
Using the definition of $\nu^{r,a}_k,$ we get 
\begin{equation}
\label{qqq}
	\widehat{\xi}_{j}^{a}(k)=\prod_{r=j+1}^{\infty}\nu^{r,a}_k=
   \prod_{r=j+1}^{\infty} \mathrm{exp}\left(-\frac{k^2+a^2}{r(r-1)a}\right)=
	\mathrm{exp}\left(-\frac{k^2+a^2}{j a}\right)=\widehat{\xi}_{j}^{a,0}(k)
\end{equation}
 for $k=-2^{j-1}+1, \dots, 2^{j-1}.$
Let us consider coefficients 
$\widehat{\xi}_{j}^{a}(k)$ and $\widehat{\xi}_{j}^{a,0}(k)$ for $|k-1/2|+1/2 \geq  2^{j-1},$ 
that is for $j\leq J-2.$  
Denoting $\nu^{r,a,0}_k:=\mathrm{exp}\left(-\frac{k^2+a^2}{r(r-1)a}\right),$ 
 recalling (\ref{xi0}), and using $|\nu^{r,a,0}_k| \leq 1,$ $|\nu^{r,a}_k| \leq 1$, we obtain
$$
\left|\widehat{\xi}_{j}^{a}(k)-\widehat{\xi}_{j}^{a,0}(k)\right|=
\left|\prod_{r=j+1}^{J-1}\nu^{r,a}_k-\prod_{r=j+1}^{J-1}\nu^{r,a,0}_k
\right|\mathrm{exp}\left(-\frac{k^2+a^2}{(J-1)a}\right)
\leq
2 \mathrm{exp}\left(-\frac{k^2+a^2}{(J-1)a}\right).
$$
Applying this estimate, (\ref{qqq}), and elementary inequalities 
$\lfloor\log_2 (k+1)\rfloor+2\leq 4 k^{1/2}$ and $k^2+a^2\geq a^{5/4}k^{3/4}$  
($a,\,k \geq 1$),  we sequentially get  
$$\|(\xi_{j}^{a})'-(\xi_{j}^{a,0})'\|^2=
\sum_{k\in \mathbb{Z}} k^2 \left|\widehat{\xi}_{j}^{a}(k)-\widehat{\xi}_{j}^{a,0}(k)\right|^2 \leq
8 \sum_{k=2^{j-1}}^{\infty} k^2 \mathrm{exp}
\left(-\frac{2(k^2+a^2)}{(\lfloor\log_2 (k+1)\rfloor+2)a}\right)\leq
 8\sum_{k=2^{j-1}}^{\infty} k^2 \mathrm{exp}\left(-\frac{1}{2}a^{1/4}k^{1/2}\right).
$$
The last expression is a remainder of a convergent series. Therefore, 
it  tends to $0$ as $j \to \infty.$ Moreover,
the last series converges uniformly on $a>1.$ So, it tends to $0$ as $a \to \infty.$
The uniformness on  $j\in \mathbb{N}$ is clear. Thus we have (\ref{juclim}).
   
3.
To check (\ref{auclim})    we use the same method as above.
 The functions $\psi_{j}^{a},$
$\eta_{j}^{a},$ $\eta_{j}^{a,0}$ (definitions are given in (\ref{uep0}), 
Remark \ref{homo}, (\ref{eta})) play the role 
of $\varphi_{j}^{a}$, $\xi_{j}^{a}$ and $\xi_{j}^{a,0}$ respectively. 
By (\ref{uep0}) and Remark \ref{homo},
$ \widehat{\eta}_{j}^{a} (k)=
\mathrm{e}^{2\pi \mathrm{i} 2^{-j-1}k}\nu^{j+1,a}_{k+2^{j}}
\widehat{\xi}_{j+1}^{a}(k).$
 Since
$\nu^{j+1,a}_{k+2^{j}}=\sqrt{1-\mathrm{exp}\left(-2(k^2+a^2)/(j(j+1)a)\right)},$ 
 $\widehat{\xi}_{j+1}^{a}(k)=\widehat{\xi}_{j+1}^{a,0}(k)=\mathrm{exp}(-(k^2+a^2)(j+1)^{-1}a^{-1})$
for $k=-2^{j-1}+1, \dots, 2^{j-1},$  
 then, recalling (\ref{eta}) we conclude (compare with (\ref{qqq}))  
$
\widehat{\eta}_{j}^{a} (k)=\widehat{\eta}_{j}^{a,0} (k)
$
as
$
k=-2^{j-1}+1, \dots, 2^{j-1}.
$
Using the same arguments as for the scaling sequence in item 2, it can be shown that
$
\lim_{j \to \infty}\|(\eta_{j}^{a})'-(\eta_{j}^{a,0})'\|=0
$
and
$
\lim_{a \to \infty}\|(\eta_{j}^{a})'-(\eta_{j}^{a,0})'\|=0
$
are fulfilled uniformly on $a>1$ and $j\in \mathbb{N}$ respectively. 
Therefore by Lemma \ref{continuous}, 
$
\lim_{j \to \infty}\sup_{a>1}|UC(\eta_{j}^{a})-UC(\eta_{j}^{a,0})|=0
$
and
$
\lim_{a \to \infty}\sup_{j>0}|UC(\eta_{j}^{a})-UC(\eta_{j}^{a,0})|=0.
$
Hence, to conclude the proof of Theorem \ref{opt} it remains to use Lemma \ref{expw}.\hfill{$\Box$}

\begin{table}[h]
\begin{center}
\begin{tabular}{c|cccccc}
$a$ & $1.1$ & $1.1$  & $1.01$ & $1.01$ & $100$ & $1000$ \\
$j$ & $10^6$ & $2\cdot 10^6$ & $5 \cdot 10^5$ &  $10^6$ & $10$ & $10$ \\
\hline
$UC(\psi^{a}_j)$ & $1.497$ & $1.498$ & $1.496$ & $1.497$ & $0.500124$ & $0.500013$ \\
\end{tabular}
\caption{Values of $UC(\psi_{j}^a)$ for particular $a$'s and $j$'s.} 
\label{tab}
\end{center}
\end{table}

\begin{figure}[thb]
		\subfigure[{}]{
		\includegraphics[width=.48\textwidth]{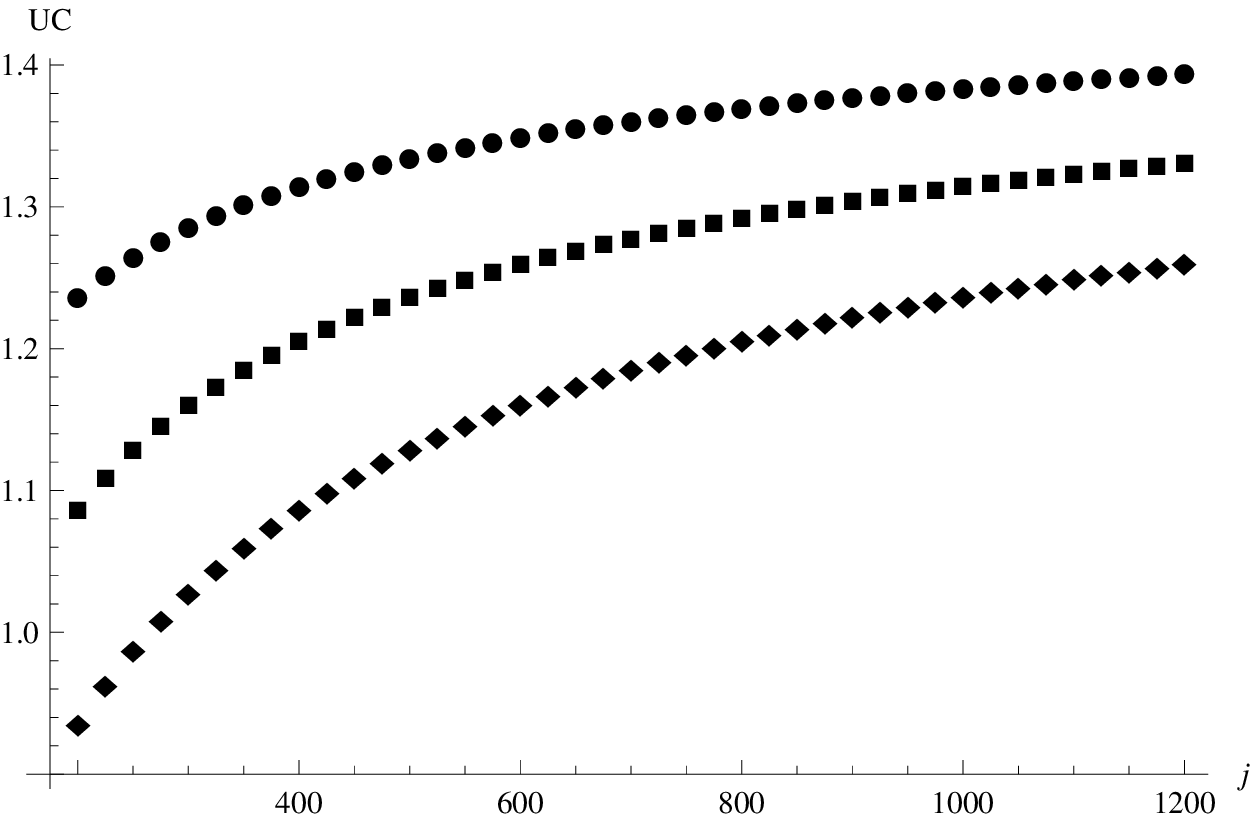}
		\label{bild1}
		}
		\subfigure[{}]{
		\includegraphics[width=.48\textwidth]{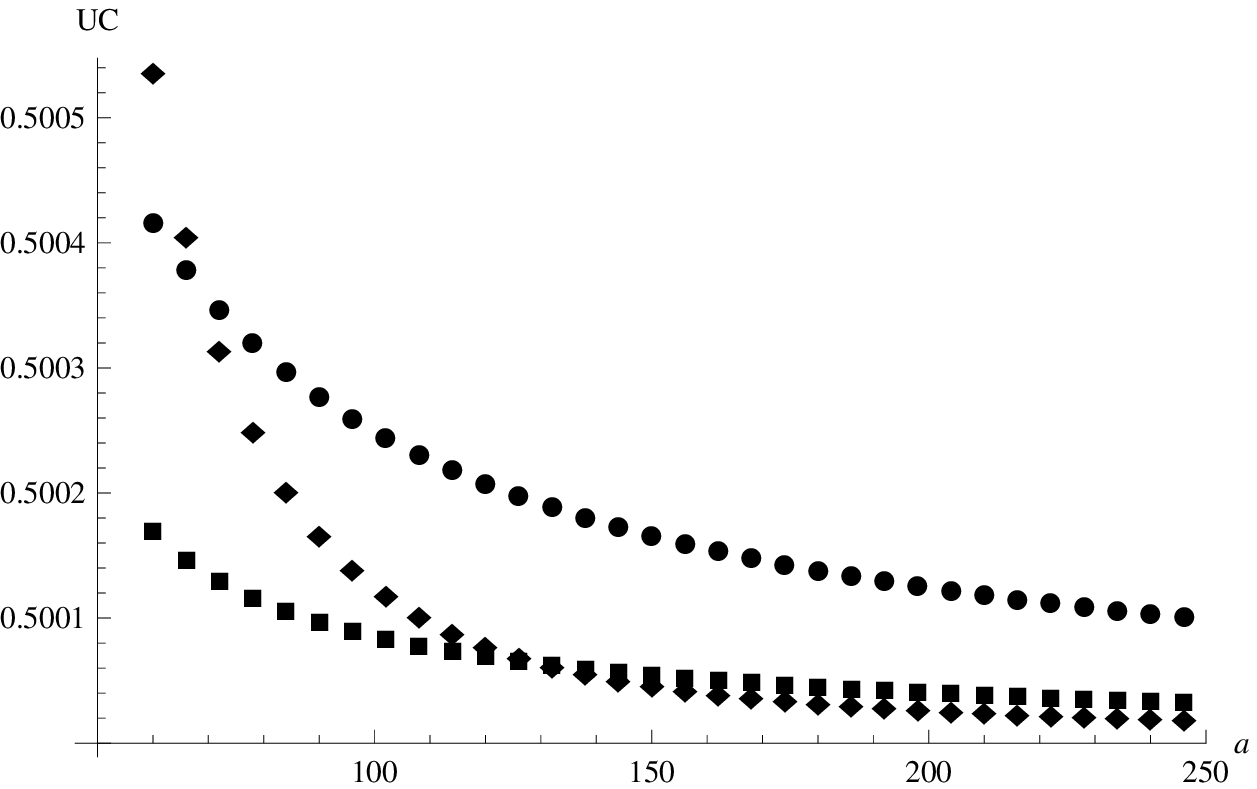}
		\label{bild2}
		}
 \caption{(a) Values of $UC(\psi_{j}^a)$ for fixed $a$: "circles", "squares", and "diamonds" correspond 
to $a=2$, $a=5$, and $a=10$ respectively. \newline 
(b) Values of $UC(\psi_{j}^a)$ for fixed $j$: "circles", "squares", and "diamonds" correspond  
to $j=5$, $j=15$, and $j=30$ respectively.}
\end{figure}

\section{Discussion}\label{disc}
In Theorem \ref{opt},  we get the optimal UC as $j \to \infty$
for the scaling sequences, but the wavelet sequences have
UCs equal to  $3/2$ only. This gives rise to a discussion. 
Let $\psi^0 \in L_2(\mathbb{R})$ be a wavelet function on the real line. 
Put 
$
\psi^p_{j,k}(x):=2^{j/2}\sum_{n \in \mathbb{Z}} \psi^0(2^j(x+n)+k).
$
The sequence $\psi^p_{j,k}$ is said to be \texttt{a periodic wavelet set generated by periodization}.  
We get the following
\begin{teo}\label{periodiz}
Suppose  $\{2^{j/2}\psi^0(2^j \cdot -k)\}_{j,k \in \mathbb{Z}}$
is a Bessel sequence and $((\psi^0)',\,\psi^0)_{L_2(\mathbb{R})}=0,$
then 
$
\lim_{j \to \infty}UC(\psi^p_{j,k}) \geq 3/2. 
$
\end{teo}

\noindent  
\textbf{Proof.}
Under aforementioned restrictions the equality $UC_H(\psi^0)\geq 3/2$  is proven 
in \cite{Bat}, \cite{Balan}. It remains to use the main result of \cite{prqurase03}, 
namely  $\lim_{j \to \infty}UC(\psi^p_{j,k})=UC_H(\psi^0).$ \hfill$\Box$ 
 
\medskip
 
These arguments motivate a conjecture: if $(\psi'_j,\,\psi_j)_{L_2(0,\,1)}=0,$
then  $\lim_{j\to \infty}UC(\psi_j)\geq 3/2$ for any periodic wavelet sequence $(\psi_j)_j$.
If this is true, the family of Parseval wavelet frames constructed in Theorem \ref{opt}
has the optimal UC. To prove the conjecture is a task for future investigation.
 
In conclusion we note that the periodization of a real line wavelet function can not 
provide a result stronger than the one in Theorem \ref{opt}. Namely, suppose 
 $f^n\in L_2(\mathbb{R}),$ $n \in \mathbb{N}$ is a sequence 
such that $\lim_{n \to \infty}UC_H(f^n)=1/2.$ 
Using the periodization we define sequences of periodic functions 
$f_{j}^{n, p}(x):=\sum_{k\in\mathbb{Z}}f^n(2^j(x+k))$ and applying results from
\cite{prqurase03} we get only
$
\lim_{j\to \infty}\lim_{n \to \infty}UC(f_{j}^{n,p})=1/2.
$    
However, it is weaker than the equalities of the form (\ref{juclim}), (\ref{auclim}).

\section*{Acknowledgments.}
The authors thank  Professor M. A. Skopina for valuable discussions.

\end{document}